\documentclass[11pt]{article}
\usepackage{amscd}
\usepackage{amsfonts}
\usepackage{amsmath}
\usepackage{amssymb}
\usepackage{amsthm}
\usepackage{bbm}
\usepackage{CJK}
\usepackage{fancyhdr}
\usepackage{graphicx}
\usepackage{hyperref}
\usepackage{indentfirst}
\usepackage{latexsym}
\usepackage{mathrsfs}
\usepackage{xypic}
\allowdisplaybreaks

\newtheorem{thm}{Theorem}[section]
\newtheorem{prop}[thm]{Proposition}
\newtheorem{defn}[thm]{Definition}

\newtheorem{lem}[thm]{Lemma}
\newtheorem{exam}[thm]{Example}

\usepackage[top=1in,bottom=1in,left=1.25in,right=1.25in]{geometry}
\def\<{\langle}
\def\>{\rangle}
\def\a{\alpha}

\def\ci{\circ}
\def\c{\cdot}

\def\D{\Delta}

\def\g{\gamma}

\def\o{\otimes}

\def\v{\varepsilon}

\def\<{\langle}
\def\>{\rangle}

\date{}
\begin{document}
\renewcommand{\baselinestretch}{1.2}
\renewcommand{\arraystretch}{1.0}
\title{\textbf{Drinfeld twists for monoidal Hom-bialgebras }}
\author{{\bf Xiaohui Zhang\footnote {E-mail: zxhhhhh@hotmail.com}
        ~and Xiaofan Zhao \footnote {Corresponding author, E-mail: zhaoxiaofan8605@gmail.com}}\\
{\small Department of Mathematics, Southeast University}\\
{\small Jiangsu Nanjing 210096, P. R. CHINA}\\}
 \maketitle

\begin{center}
\begin{minipage}{14.cm}
\noindent\textbf {\bf Abstract.}
The aim of this paper is to define and study Drinfeld twists for monoidal Hom-bialgebras. We show that a new Hom-bialgebra could be constructed by
changing the coproduct of a monoidal Hom-bialgebra via a Drinfeld twist, and this construction preserves $R$-matrixes if there exist one. Moreover, their representation categories are monoidal isomorphic.
\\


\end{minipage}
\end{center}
\section{\textbf{Introduction}}

\setcounter{equation}{1}

In 2006, Hartwig, Larsson and Silvestrov introduced the Hom-Lie algebras when they concerned about the $q$-deformations of Witt and Virasoro algebras (see \cite{JDS}). In a Hom-Lie algebra, the Jacobi identity is replaced by the so called Hom-Jacobi identity via a homomorphism. Hom-associative algebras, the corresponding structure of associative algebras, were introduced by Makhlouf and Silvestrov in \cite{AS2}. The associativity of a Hom-algebra is twisted by an endomorphism (here we call it the Hom-structure map). The
generalized notions, Hom-bialgebras, Hom-Hopf algebras were developed in \cite{AS1}, \cite{AS3}, \cite{AS4}.
Further research on various Hom-Lie structures and Hom-type algebras by many schlors could be found in \cite{AF}, \cite{dy3}, \cite{MPJ}, \cite{hnh}. Quasitriangular Hom-bialgebras were considered by Yau (\cite{dy2}, \cite{dy3}), which provided a solution of the quantum Hom-Yang-Baxter euqation, a twisted version of the quantum Yang-Baxter equation (\cite{dy4}, \cite{dy5}).

The notions of Hom-categories and monoidal Hom-Hopf algebras were introduced by Caenepeel and Goyvaerts (\cite{CG}) in order to provide a categorical approach to Hom-type algebras. In a Hom-category, the associativity and unit constraints are twisted by the Hom-structure maps. A (co)monoid in the Hom-category is a Hom-(co)algebra, and a bimonoid in the Hom-category is a monoidal Hom-bialgebra. Further research on monoidal Hom-bialgebras can be found in \cite{cy1}, \cite{cy2}, and \cite{ls}.

Moreover, through a direct computation (see Example \ref{bm} and Example \ref{bh} for details), we can get that there is a one to one correspondence between the collection the monoidal Hom-bialgebras over a commutative ring $k$, and the collection of the unital Hom-bialgebra over $k$ which Hom-structure map is a bijection. Does there another way to get a Hom-bialgebra through a given monoidal Hom-bialgebra? Is there any relation between their representation categories? This is the motivation of the present article. In order to investigate these questions, we
introduce the definition of Drinfeld twists for monoidal Hom-bialgebras, and construct a new Hom-bialgebra via a Drinfeld twist.

A Drinfeld twist for a Hopf algebra $H$ is an invertible element $\sigma \in H \o H$, satisfying the 2-cocycle condition
$$
(\sigma \o 1)(\Delta \o id)(\sigma) = (1\o \sigma)(id \o \Delta)(\sigma).
$$
Note that our definition of a Drinfeld twist is inverse with Drinfeld's (see \cite{dr} and \cite{eg}).
In our paper, we always assume that $\sigma$ is normalized,
i.e. $(\varepsilon \o id)(\sigma)=( id \o\varepsilon)(\sigma)=1_H$.

The twisting elements or twists were first introduced by Drinfeld \cite{dr} on quasi-Hopf algebras, in
order to twist the coproduct without changing its product. They have become an important
tool in the classification of finite-dimensional Hopf algebras (\cite{eg}). The twisting elements for a Hom-bialgebra have been discussed in \cite{mld}.

The paper is organized as follows. In Section 2 we recall some notions of monoidal Hom-type algebras and Hom-type algebras. In section 3, we describe the category of representations of a monoidal Hom-bialgebra which is more generalized than Caenepeel's, and give the definition of quasitriangular monoidal Hom-bialgebras. In section 4, we introduce the notion of the Drinfeld twists for a monoidal Hom-bialgebra, and construct a Hom-bialgebra by changing the coproduct via a Drinfeld twist. If the monoidal Hom-bialgebra is quasitriangular, then the Hom-bialgebra which we obtained is also quasitriangular.   Furthermore, we show that their representation categories are monoidal isomorphic.

\section{\textbf{Preliminaries}}

Throughout the paper, we let $k$ be a fixed
commutative ring and $char(k)=0$. All algebras are supposed to be over $k$. For the comultiplication
 $\D $ of a $k$-module $C$, we use the Sweedler-Heyneman's notation:
 $$
 \Delta(c)=c_{1}\o c_{2},
 $$
for any $c\in C$. $\tau$ means the flip map $\tau(a \o b) = b \o a$.
When we say a "Hom-algebra" or a "Hom-coalgebra", we mean the unital Hom-algebra and counital Hom-coalgebra.

\vskip 0.5cm
 {\bf 2.1. Monoidal Hom-bialgebras and monoidal Hom-Hopf algebras.}
\vskip 0.5cm

In this section, we will review several definitions and notations related to monoidal Hom-bialgebras (see \cite{CG}).

Let $\mathcal{C}$ be a category. We introduce a new category
$\widetilde{\mathscr{H}}(\mathcal{C})$ as follows: the objects are
couples $(M, \a_M)$, with $M \in \mathcal{C}$ and $\a_M \in
Aut_{\mathcal{C}}(M)$. A morphism $f: (M, \a_M)\rightarrow (N, \a_N)$
is a morphism $f : M\rightarrow N$ in $\mathcal{C}$ such that $\a_N
\circ f= f \circ \a_M$.

Specially, let $\mathscr{M}_k$ denote the category of $k$-modules.
~$\mathscr{H}(\mathscr{M}_k)$ will be called the Hom-category
associated to $\mathscr{M}_k$. If $(M,\a_M) \in \mathscr{M}_k$, then
$\a_M: M\rightarrow M$ is obviously an isomorphism in
~$\mathscr{H}(\mathscr{M}_k)$. It is easy to show that
~$\widetilde{\mathscr{H}}(\mathscr{M}_k)$ =
(~$\mathscr{H}(\mathscr{M}_k),~\otimes,~(k, id),~\widetilde{a},
~\widetilde{l},~\widetilde{r}))$ is a monoidal category by
Proposition 1.1 in \cite{CG}:

$\bullet$ the tensor product of $(M,\a_M)$ and $(N,
\a_N)$ in ~$\widetilde{\mathscr{H}}(\mathscr{M}_k)$ is given by the
formula $(M, \a_M)\otimes (N, \a_N) = (M\otimes N, \a_M \otimes \a_N)$.

$\bullet$ for any $m \in M$, $n \in N$, $p \in P$, the associativity is given by the formulas
$$
\widetilde{a}_{M,N,P}((m\o n)\o p)=\a_M(m)\o (n\o \a_P^{-1}(p)),$$

$\bullet$ for any $m \in M$, the unit
constraints are given by the formulas
$$
\widetilde{l}_{M}(x\o m)=\widetilde{r}_{M}(m\o x)=x\a_M(m).
$$

A \emph{monoidal Hom-algebra} over $k$  is an
object $(A, \alpha)\in
\widetilde{\mathscr{H}}(\mathscr{M}_k)$ together with a $k$-linear
map $m_A: A\o A\rightarrow A$ and an element $1_A\in A$ such that for all $a,b,c\in A$,
\\
$\left\{\begin{array}{l}
(1)~~\a(ab) = \a(a)\a(b);\\
(2)~~\a(a)(bc) = (ab)\a(c);\\
(3)~~\a(1_A) = 1_A;\\
(4)~~1_Aa = a1_A = \a(a).
\end{array}\right.$
\\

\emph{A morphism $f:A\rightarrow B$ of monoidal Hom-algebras} is a linear map such that $\a_B\circ f = f\circ \a_A$, $f(1_A) = 1_B$ and $\mu_B\circ (f \o f) = f\circ \mu_A$.

 A \emph{monoidal Hom-coalgebra} over $k$ is
an object $(C,\g)\in \widetilde{\mathscr{H}}(\mathscr{M}_k)$
together with $k$-linear maps $\Delta:C\rightarrow C\otimes C,~~~
\D(c)=c_{1}\o c_{2}$ and
$\a:C\rightarrow C$ such that for all $c\in C$,
\\
$\left\{\begin{array}{l}
(1)~~\Delta(\g(c)) = \g(c_1)\g(c_2);\\
(2)~~\g^{-1}(c_1) \o \Delta(c_{2}) = \Delta(c_1) \o \g^{-1}(c_2);\\
(3)~~\varepsilon \circ \g = \varepsilon;\\
(4)~~\varepsilon(c_1)c_2 = c_1\varepsilon(c_2) = \g^{-1}(c).
\end{array}\right.$
\\

\emph{A morphism $f:C\rightarrow D$ of monoidal Hom-coalgebras} is a linear map such that $\g_D\circ f = f\circ \g_C$, $\varepsilon_C = \varepsilon_D\circ f$ and $\Delta_D\circ f = (f \o f)\circ \Delta_C$.

 A \emph{monoidal Hom-bialgebra}
$H=(H,\a,m,\eta, \Delta,\varepsilon)$  is a bimonoid in
$\widetilde{\mathscr{H}}(\mathscr{M}_k)$. This means that $(H, \a,
m,\eta)$ is a monoidal Hom-algebra, $(H,\alpha,\Delta,\varepsilon)$ is a monoidal Hom-coalgebra
and $\D$, $\v$ are morphisms of monoidal Hom-algebras preserving unit.

\begin{exam}\label{bm}
Suppose $(A,m,\eta,\Delta,\varepsilon)$ is a bialgebra over $k$ endowed with a bialgebra isomorphism $\a:A\rightarrow A$. Then $(A,\a,\a\circ m,\eta,\Delta\circ \a^{-1},\varepsilon)$ is a monoidal Hom-bialgebra over $k$. We denote this monoidal Hom-bialgebra by ${}^\a A$.

Conversely, if $(H,\a,m,\eta,\Delta,\varepsilon)$ is a monoidal Hom-bialgebra, then $(H,\a^{-1}\circ m,\eta,\Delta\circ \a,\varepsilon)$ is a bialgebra over $k$. We write this bialgebra for ${}_\a H$.

Thus we immediately get a bijective map $A\rightarrow {}^\a A$ between
the collection of all bialgebras over $k$ endowed with an invertible endomorphism on it and the collection of all monoidal Hom-bialgebras over $k$.
\end{exam}

 A \emph{monoidal Hom-Hopf algebra} over $k$ is
a monoidal Hom-bialgebra $(H, \alpha)$ together with a linear map
$S:H\rightarrow H$ in $\widetilde{\mathscr{H}}(\mathscr{M}_k)$ (called \emph{the antipode}) such
that
$$S\ast id=id\ast S=\eta\varepsilon,~S\circ\alpha=\alpha\circ S.$$

Recall from (Proposition 2.9, \cite{CG}) that for all $a, b \in H$, $S$ satisfies
$$
S(ab) = S(b)S(a),~~S(1_H) = 1_H,
$$
$$
\Delta(S(a)) = S(a_2) \o S(a_1),~~\varepsilon\circ S=\varepsilon.
$$

Note that the map $A\rightarrow {}^\a A$ given in Example \ref{bm} is also a bijection between the collection of $k$-Hopf algebra $(A,m,\eta,\Delta,\varepsilon,S)$ which endowed with an invertible endomorphism and the collection of $k$-monoidal Hom-Hopf algebra
 ${}^\a A=(A,\a,\a \ci m, \eta,\Delta \ci \a^{-1},\varepsilon,S)$.

\vskip 0.5cm
 {\bf 2.2. Hom-bialgebras and Hom-Hopf algebras.}
\vskip 0.5cm

In this section, we will review several definitions and notations related to Hom-bialgebras.

Recall from \cite{AS2} that
a \emph{Hom-algebra} over $k$ is a quadruple $(A,\mu,\eta,\alpha)$, in which $A$ is a $k$-module, $\alpha:A\rightarrow A$, $\mu: A \o A\rightarrow A$ and $\eta:k \rightarrow A$ are linear maps, with notation $\mu(a \o b) = ab$ and $\eta(1_k) = 1_A$, satisfying the following conditions, for all $a,b,c \in A$:\\
$\left\{\begin{array}{l}
(1)~~\a(ab) = \a(a)\a(b);\\
(2)~~\a(a)(bc) = (ab)\a(c);\\
(3)~~\a(1_A) = 1_A;\\
(4)~~1_Aa = a1_A = \a(a).
\end{array}\right.$
\\

Note that a monoidal Hom-algebra is also a Hom-algebra. Conversely, a Hom-algebra is a monoidal Hom-algebra if its Hom-structure map is invertible.

\emph{A morphism $f:A\rightarrow B$ of Hom-algebras} is a linear map such that $\a_B\circ f = f\circ \a_A$, $f(1_A) = 1_B$ and $\mu_B\circ (f \o f) = f\circ \mu_A$.

Recall from \cite{AS1} that
a \emph{Hom-coalgebra} over $k$ is a quadruple $(C,\Delta,\varepsilon,\a)$, in which $C$ is a $k$-module, $\a:C\rightarrow C$, $\Delta: C\rightarrow C\o C$ and $\varepsilon:C\rightarrow k$ are linear maps, with notation $\Delta(c) = c_1 \o c_2$, satisfying the following conditions for all $c \in C$:\\
$\left\{\begin{array}{l}
(1)~~\Delta(\a(c)) = \a(c_1)\o \a(c_2);\\
(2)~~\a(c_1) \o \Delta(c_{2}) = \Delta(c_1) \o \a(c_2);\\
(3)~~\varepsilon \circ \a = \varepsilon;\\
(4)~~\varepsilon(c_1)c_2 = c_1\varepsilon(c_2) = \a(c).
\end{array}\right.$
\\

\emph{A morphism $f:C\rightarrow D$ of Hom-coalgebras} is a linear map such that $\a_D\circ f = f\circ \a_C$, $\varepsilon_C = \varepsilon_D\circ f$ and $\Delta_D\circ f = (f \o f)\circ \Delta_C$.

Recall from \cite{AS4} that
a \emph{Hom-bialgebra}
$B=(B,\a,m,\eta, \Delta,\varepsilon)$ over $k$ is a sextuple $(B,\mu,\eta$, $\Delta$, $\varepsilon$, $\a)$, in which $(B, \a,
m,\eta)$ is a Hom-algebra, $(B,\alpha,\Delta,\varepsilon)$ is a Hom-coalgebra,
and $\D$, $\v$ are morphisms of Hom-algebras preserving unit.

\begin{exam}\label{bh}
Suppose $(A,m,\eta,\Delta,\varepsilon)$ is a bialgebra over $k$ endowed with a bialgebra map $\a:A\rightarrow A$. Then $(A,\a,\a\circ m,\eta,\Delta\circ \a,\varepsilon)$ is a Hom-bialgebra over $k$. We denote this Hom-bialgebra by $A^\a$.

Conversely, if $(B,\a,m,\eta,\Delta,\varepsilon)$ is a Hom-bialgebra and $\a$ is invertible, then $(B,\a^{-1}\circ m,\eta,\Delta\circ \a^{-1},\varepsilon)$ is a bialgebra over $k$. We write this bialgebra for $B_\a$.

Thus we immediately get a bijective map $A\rightarrow A^\a$ between
the collection of all bialgebras over $k$ endowed with an invertible endomorphism on it, and the collection of all Hom-bialgebras with invertible Hom-structure maps.
\end{exam}

\begin{defn}
 A \emph{Hom-Hopf algebra} over $k$ is
a Hom-bialgebra $(H, \alpha)$ together with a $k$-linear map
$S:H\rightarrow H$ (called \emph{the antipode}) such
that
$$S\ast id=id\ast S=\eta\varepsilon,~S\circ\alpha=\alpha\circ S.$$
\end{defn}

\begin{lem}
Let $(H, \a)$ be a Hom-Hopf algebra. If $\a$ is invertible, then for all $a, b \in H$, $S$ satisfies
$$
S(ab) = S(b)S(a),~~S(1_H) = 1_H,
$$
$$
\Delta(S(a)) = S(a_2) \o S(a_1),~~\varepsilon\circ S=\varepsilon.
$$
\end{lem}

{\bf Proof.}
We will only prove the first statement. We compute as follows.
$$\aligned
S(ab)&=S(\a^{-1}(a_1)\a^{-1}(b_1))\varepsilon(b_2)\varepsilon(a_2)\\
&=S(\a^{-2}(a_1)\a^{-2}(b_1))\varepsilon(b_2)(\a^{-3}(a_{21})S(\a^{-3}(a_{22})))\\
&=S(\a^{-2}(a_1)\a^{-2}(b_1))( ( \a^{-4}(a_{21}) (\a^{-5}(b_{21})S(\a^{-5}(b_{22})))) S(\a^{-3}(a_{22})))\\
&=S(\a^{-2}(a_1)\a^{-2}(b_1))( (\a^{-4}(a_{21}) \a^{-4}(b_{21})) (S(\a^{-4}(b_{22})) S(\a^{-4}(a_{22}))) )\\
&=( S(\a^{-4}(a_{11}) \a^{-4}(b_{11})) (\a^{-4}(a_{12}) \a^{-4}(b_{12})) ) (S(\a^{-2}(b_{2})) S(\a^{-2}(a_{2})))\\
&=\varepsilon(a_1)\varepsilon(b_1)S(\a^{-1}(b_{2})) S(\a^{-1}(a_{2}))=S(b)S(a).
\endaligned$$
 $\hfill \Box$

Note that the map $A\rightarrow A^\a$ given in Example \ref{bh} is also a bijection between the collection of all $k$-Hopf algebras endowed with an invertible endomorphism and the collection of all $k$-Hom-Hopf algebras with invertible Hom-structure maps.

\begin{defn}
Let $(B,\a)$ be a Hom-bialgebra. If there exists an invertible element $R \in B \o B$, such that the following conditions hold:\\
$\left\{\begin{array}{l}
(q1)~~(\a \o \a)R =R;\\
(q2)~~R\Delta(x) = \Delta^{op}(x)R;\\
(q3)~~R^{(1)}_1 \o R^{(1)}_2 \o \a(R^{(2)}) = \a(r^{(1)}) \o \a(R^{(1)}) \o r^{(2)}R^{(2)};\\
(q4)~~\a(R^{(1)}) \o R^{(2)}_1 \o R^{(2)}_2 = r^{(1)}R^{(1)} \o \a(R^{(2)}) \o \a(r^{(2)}),
\end{array}\right.$
\\
for any $x \in B$, where $r=R = R^{(1)} \o R^{(2)} = r^{(1)} \o r^{(2)}$, then $R$ is called an \emph{$R$-matrix} of $B$, $(B,\a,R)$ is called a \emph{quasitriangular Hom-bialgebra}.
\end{defn}

Note that the above definition of quasitriangular Hom-bialgebras is slightly different from Yau's (see \cite{dy2}). In order to make sure the representation category of a quasitriangular Hom-bialgebra is braided, we need $R$ is invertible and Eq.(q1) is hold.

 Let $(B,\alpha)$ be a Hom-algebra.  A left \emph{$(B,\alpha)$-Hom-module} is a triple $(M,\a_M, \theta_M)$, where $M$ is a $k$-module, $\a_M:M\rightarrow M$ and $\theta_M:B \o M\rightarrow M$ are linear maps with notation $\theta_M(b \o m) = b\cdot m$, satisfying the following conditions, for all $b,b'\in B$, $m\in M$
\\
$\left\{\begin{array}{l}
(1)~~\a(b\cdot m) = \a(b)\cdot\a_M(m);\\
(2)~~\a(b)\cdot(b'\cdot m) = (bb')\cdot\a_M(m);\\
(3)~~1_B\cdot m = \a_M(m).
\end{array}\right.$
\\

A \emph{morphism $f:M\rightarrow N$ of $B$-modules} is a $k$-linear map such that $\a_N\circ f = f\circ \a_M$ and $\theta_N\circ (id_A \o f) = f\circ \theta_M$.

\section{\textbf{The representations of monoidal Hom-bialgebras}}

\setcounter{equation}{0}

Let $(H, \a)$ be a monoidal Hom-bialgebra. A category $Rep^{i,j}(H)$ is defined as follows for any fixed $i,j \in \mathbb{Z}$ (the domain of integrals):

$\bullet $ $Rep^{i,j}(H)$ is the category of left Hom-modules of the monoidal Hom-algebra $H$ and the morphisms of $H$-modules;

$\bullet $ the tensor product $M \o N$ for $(M, \a_M), (N,\a_N) \in Rep^{i,j}(H)$ is obtained by
$(M \o N, \a_M \o \a_N)$
with the action of $H$ given by
$$h\cdot (m \o n) = \a^{i}(h_1)\cdot m \o \a^j(h_2)\cdot n, $$
where $i,j \in \mathbb{Z}$, $h \in H $, $m \in M$, $n \in N$;

$\bullet $ the tensor product of two arrows $f,g \in Rep^{i,j}(H)$ is given by the tensor product of $k$-linear morphisms, i.e. the forgetful functor from $Rep^{i,j}(H)$ to the category of $k$-modules is faithful;

$\bullet $ for any $\lambda \in k$, $(k,id_k)$ is an object in $Rep^{i,j}(H)$ with the action
$$
h\cdot \lambda = \varepsilon(h)\lambda.
$$

\begin{lem}
$k$ is the unit object of the tensor product $\o$ in $Rep^{i,j}(H)$.
\end{lem}

{\bf Proof.}
Firstly, it is easy to check that $(k, id_k)$ and $(M \o N, \a_M \o \a_N)$ under the $H$-action defined above are objects in $Rep^{i,j}(H)$.

Secondly, for any $M \in Rep^{i,j}(H)$, define a $k$-linear map
$$
l_{M}:~k \o M \rightarrow M,~~~~\lambda \o m\mapsto \lambda\a_M^{-j+1}(m),~~\lambda \in k,~m \in M.
$$
Obviously $l$ is natural, and we have $\a_M \circ l_M = l_M \circ(id_k \o \a_M)$. Furthermore, $l_M$ is $H$-linear, actually,
$$\aligned
l_{M}(h\cdot(\lambda \o m))= & l_{M}(\varepsilon (h_{1})\lambda \o \a^j(h_{2})\cdot m) \\
&=\lambda\a^{-j+1}_M(\a^{j-1}(h)\cdot m) \\
&=\lambda h \cdot \a_M^{-j+1}(m) = h\cdot l_{M}(\lambda \o m).
\endaligned$$
The inverse of $l$ is given by
$$
l_{M}^{-1}:~M\rightarrow k \o M,~~~~ m\mapsto 1_k \o \a_M^{j-1}(m),~~m \in M.
$$

Similarly, we define the $k$-linear maps
$$
r_{M}:~M \o k  \rightarrow M,~~~~ m \o \lambda \mapsto \lambda \a_M^{-i+1}(m),~~\lambda \in k ,~m \in M,
$$
and
$$
r_{M}^{-1}:~M\rightarrow M \o k ,~~~~ m\mapsto \a^{i-1}_M(m) \o 1_k ,~~m \in M.
$$
It is easy to check that $r$ is a natural isomorphism with the inverse $r^{-1}$.

This completes the proof.
 $\hfill \Box$

\begin{thm}
$Rep^{i,j}(H)$ is a monoidal category.
\end{thm}

{\bf Proof.}
Firstly, for any $M,N,P,Q \in Rep^{i,j}(H)$, define an associativity constraint by
$$
a_{M,N,P}((m \o n) \o p) = \a^{-i+1}_M(m) \o(n \o \a^{j-1}_P(p)),~m \in M,~n \in N,~p \in P.
$$
Obviously that $a$ is natural and satisfies $a_{M,N,P}\circ(\a_M \o (\a_N \o \a_P)) = ((\a_M \o \a_N)\o \a_P)\circ a_{M,N,P}$. For any $h \in H$, we obtain
$$\aligned
&~~~~a_{M,N,P}(h\cdot ((m \o n) \o p))\\
&= \a_M^{-i+1}(\a^{2i}(h_{11})\cdot m )\o(\a^{i+j}(h_{12})\cdot n \o \a_P^{j-1}(\a^{j}(h_{2})\cdot p))\\
&=\a^{i}(h_1)\cdot \a_M^{-i+1}(m) \o ( \a^{i+j}(h_{21})\cdot n \o \a^{2j}(h_{22})\cdot \a_P^{j-1}(p))\\
&=h\cdot(\a_M^{-i+1}(m) \o (n \o \a_P^{j-1}(p) ))\\
&= h\cdot (a_{M,N,P}((m \o n) \o p) ),
\endaligned$$
thus $a_{M,N,P}$ is $H$-linear. Since $a$ is invertible, $a$ is a natural isomorphism in $Rep^{i,j}(H)$.

Secondly, one can see that $a$ satisfies the Pentagon Axiom. Actually,
$$\aligned
&~~~~((id_M \o a_{N,P,Q})\circ a_{M,N\o P, Q} \circ (a_{M,N,P} \o id_Q))(((m \o n) \o p) \o q)\\
&=(id_M \o a_{N,P,Q})(\a_M^{-2i+2}(m) \o ((n \o \a_P^{j-1}(p)) \o \a_Q^{j-1}(q)))\\
&=\a_M^{-2i+2}(m) \o (\a_N^{-i+1}(n) \o (\a_P^{j-1}(p) \o \a_Q^{2j-2}(q)))\\
&=a_{M,N,P \o Q}((\a_M^{-i+1}(m) \o \a_N^{-i+1}(n)) \o (p \o \a_Q^{j-1}(q)))\\
&=(a_{M \o N,P,Q}\circ a_{M,N,P \o Q})(((m \o n) \o p) \o q).
\endaligned$$

At last, it is also a direct check to prove that $a,l,r$ satisfy the Triangle Axiom.
This completes the proof.
$\hfill \Box$

\begin{defn}
Let $(H,\a)$ be a monoidal Hom-bialgebra. If there exists an invertible element $R = R^{(1)} \o R^{(2)} \in H \o H$, such that the following conditions hold:\\
$\left\{\begin{array}{l}
(Q1)~~(\a \o \a)R =R;\\
(Q2)~~R\Delta(h) = \Delta^{op}(h)R;\\
(Q3)~~R^{(1)}_1 \o R^{(1)}_2 \o R^{(2)} = r^{(1)} \o R^{(1)} \o r^{(2)}R^{(2)};\\
(Q4)~~R^{(1)} \o R^{(2)}_1 \o R^{(2)}_2 = r^{(1)}R^{(1)} \o R^{(2)} \o r^{(2)},
\end{array}\right.$
\\
where $h \in H$, $r=R = R^{(1)} \o R^{(2)} = r^{(1)} \o r^{(2)}$, then $R$ is called an \emph{$R$-matrix} of $H$. $(H,\a,R)$ is called a \emph{quasitriangular monoidal Hom-bialgebra}.
\end{defn}

\begin{exam}
Let $(A,R)$ be a quasitriangular bialgebra over $k$ and $\a:A\rightarrow A$ be an invertible bialgebra homomorphism. If $R$ satisfies $(\a \o \a)(R)=R$,
then $({}^\a A,\a,R)$ in Example \ref{bm} is a quasitriangular monoidal Hom-bialgebra.

Conversely, if $(H,\a,R)$ is a quasitriangular monoidal Hom-bialgebra, then $({}_\a H, R)$ in Example \ref{bm} is a quasitriangular bialgebra over $k$.
\end{exam}

%

\begin{prop}
Let $(H,\a,R)$ be a quasitriangular monoidal Hom-bialgebra. Then $R$ satisfies the quantum Hom-Yang-Baxter equations
$$(R_{12}R_{13})R_{23} = R_{23}(R_{13}R_{12}),~~~R_{12}(R_{13}R_{23}) =(R_{23}R_{13})R_{12},$$
where $R_{12} = R \o 1_H$, $R_{23} = 1_H \o R$, $R_{13}=(\tau \o id)R_{23}$.
\end{prop}

{\bf Proof.}  Straightforward.
$\hfill \Box$

\begin{thm}
 Let $(H, \a)$ be a monoidal Hom-bialgebra. For the fixed elements $R,S \in H \o H$, define maps
$$c_{M,N}:M \o M\rightarrow N \o M,~~m \o n\mapsto \a^i(R^{(2)})\cdot \a_N^{i-j-1}(n) \o \a^j(R^{(1)})\cdot \a_M^{j-i-1}(m),$$
and
$$c'_{M,N}:N \o M\rightarrow M \o M,~~n \o m\mapsto \a^i(S^{(1)})\cdot \a_M^{i-j-1}(m) \o \a^j(S^{(2)})\cdot \a_N^{j-i-1}(n),$$
for any $(M,\a_M),(N,\a_N) \in Rep^{i,j}(H)$, then $c$ is a braiding in $Rep^{i,j}(H)$ with the inverse $c'$ if and only if $R$ is an $R$-matrix with the inverse $S$.
\end{thm}

{\bf Proof.}
It is easy to check that $c$ and $c'$ are natural in $Rep^{i,j}(H)$.\\

(1). If $R$ satisfies $(\a \o \a)R =R$, then we immediately get that
$$\aligned
c_{M,N}\ci \a_{M,N}=\a_{N \o M}\ci c_{M,N}.
\endaligned$$
If $R$ also satisfies $R\Delta(h) = \Delta^{op}(h)R$ for any $h \in H$, then for any $m \in M$, $n \in N$, we have
$$\aligned
c_{M,N}(h\cdot(m \o n)) &= \a^i(R^{(2)}) \c \a_N^{i-j-1}(\a^{j}(h_2)\c n) \o \a^j(R^{(1)}) \c \a_M^{j-i-1}(\a^{i}(h_1)\c m)\\
&= (\a^{i-1}(R^{(2)})\a^{i-1}(h_2))\cdot a_N^{i-j}(n) \o (\a^{j-1}(R^{(1)})\a^{j-1}(h_1))\cdot a_M^{j-i}(m) \\
&=\a^{i}(h_1)\c(\a^{i}R^{(2)})\cdot a_N^{i-j-1}(n) \o  \a^{j}(h_2)\c(\a^{j}R^{(1)})\cdot a_M^{j-i-1}(m)\\
& = h\cdot(c_{M,N}(m \o n)),
\endaligned$$
thus $c_{M,N} \in Mor(Rep^{i,j}(H))$.

Conversely, if $c_{M,N}\in Mor(Rep^{i,j}(H))$, then $c$ is $H$-linear and satisfies $c_{M,N}\circ \a_{M\o N} = \a_{N \o M} \circ c_{M,N}$. Take $M=N=H$ and $m=n=1_H$, then we directly get Eq.(Q1) and Eq.(Q2) for any $h \in H$.\\

(2). If $S$ is the inverse of $R$, thus $S$ satisfies $(\a \o \a)(S)=S$. Then for any $m \in M$, $n \in N$, we obtain
$$\aligned
&~~~~c'_{M,N}(c_{M,N}(m \o n)) \\
&= \a^{i}(S^{(1)})\cdot \a^{i-j-1}_M(\a^j(R^{(1)})\cdot \a_M^{j-i-1}(m)) \o \a^{j}(S^{(2)})\cdot \a^{j-i-1}_N(\a^i(R^{(2)})\cdot \a_N^{i-j-1}(n))\\
&= \a^{i-1}(S^{(1)}R^{(1)}) \c \a_M^{-1}(m) \o\a^{j-1}(S^{(2)}R^{(2)}) \c \a_N^{-1}(n)\\
&= m \o n.
\endaligned$$
Similarly, we have
$$\aligned
c_{M,N}\ci c'_{M,N}=id_{N \o M},
\endaligned$$
hence $c'$ is the inverse of $c$.

Conversely, if $c'$ is the inverse of $c$, take $M=N=H$ and $m=n=1_H$, then we immediately get $RS=SR = 1_H \o 1_H$.\\

(3). Assume that $R$ satisfies Eq.(Q4), thus we obtain
$$\aligned
&~~~~(a_{N,P,M}\circ c_{M,N \o P}\circ a_{M,N,P})((m \o n) \o p)\\
&=a_{N,P,M}((\a^{2i}(R^{(2)}_1)\cdot\a_N^{i-j-1}(n) \o \a^{i+j}(R^{(2)}_2)\cdot \a_P^{i-2}(p)) \o \a^j(R^{(1)})\cdot \a_M^{j-2i}(m))\\
&=\a^{i+1}(R^{(2)})\cdot\a_N^{-j}(n)\o(\a^{i+j}(r^{(2)})\cdot \a_P^{i-2}(p) \o \a^{2j-1}(r^{(1)}R^{(1)}) \cdot \a_M^{2j-2i-1}(m)))\\
&=\a(R^{(2)})\cdot\a_N^{-j}(n)\o(\a^{i}(r^{(2)})\cdot \a_P^{i-2}(p) \o \a^{j}(r^{(1)})\c(\a^{2j-i-1}(R^{(1)}) \cdot \a_M^{2j-2i-2}(m)))\\
&=(id_N \o (c_{M,P}))(\a(R^{(2)})\cdot\a_N^{-j}(n)\o (\a^{j}(R^{(1)})\cdot \a_M^{j-i-1}(m) \o \a_P^{j-1}(p)))\\
&=((id_N \o c_{M,P})\circ a_{N,M,P}\circ (c_{M,N} \o id_P))((m \o n) \o p),
\endaligned$$
for any $M,N,P \in Rep^{i,j}(H)$.

Similarly, if $R$ satisfies Eq.(Q3), then $c$ satisfies
$$a^{-1}_{P,M,N}\circ c_{M \o N,P} \circ a^{-1}_{M,N,P} = (c_{M,P} \o  id_N)\circ a^{-1}_{M,P,N}\circ(id_M \o c_{N,P}).$$

Conversely, if $c$ is a braiding, then take $M=N=P=H$ and $m=n=p=1_H$, it is a direct computation to get Eq.(Q3) and Eq.(Q4).
\\

Combining (1)-(3), the conclusion holds.
$\hfill \Box$

\begin{thm}
 Let $(H, \a)$ be a monoidal Hom-bialgebra. For any $i,j,i',j' \in \mathbb{Z}$, $Rep^{i,j}(H)$ and $Rep^{i',j'}(H)$ are monoidal isomorphic.
 Moreover, if $H$ is a quasitriangular monoidal Hom-bialgebra, then $Rep^{i,j}(H)$ and $Rep^{i',j'}(H)$ are isomorphic as braided categories.
\end{thm}

{\bf Proof.}
Define a functor $$F=(F,F_2,F_0):(Rep^{i,j}(H)\o, k, a, l, r)\rightarrow (Rep^{i'j'}(H),\o',k,a',l',r')$$
by
 $$F(M): = M,~~ F(f): = f, ~~F_0=id_{k},
$$
$$
F_2(M,N):F(M) \o' F(N) \rightarrow F(M \o N),~~m \o' n\mapsto \a^{i-i'}_M(m) \o \a^{j-j'}_N(n),
$$
where $M,N \in Rep^{i,j}(H)$, $f:M\rightarrow N  \in Mor(Rep^{i,j}(H))$, and $\o',a',l',r'$ mean the corresponding structures in $Rep^{i',j'}(H)$.

Obviously $F_2$ is natural and compatible with the Hom structure map.

Firstly, since $F_2(M,N)$ is invertible, and
$$\aligned
F_2(M,N)(h\cdot(m \o' n))&=F_2(M,N)(\a^{i'}(h_1)\cdot m \o' \a^{j'}(h_2)\cdot n)\\
&=h\cdot(\a^{i-i'}_M \o \a^{j-j'}_N(n))= h\cdot F_2(M,N)(m \o' n),
\endaligned$$
$F_2$ is a natural isomorphism in $Rep^{i',j'}(H)$.

Secondly, we have
$$\aligned
&~~~~(F_2(M,N \o P)\circ(id_{F(M)} \o' F_2(N,P))\circ(a_{F(M),F(N),F(P)}))((m\o' n)\o' p)\\
&=(F_2(M,N \o P)\circ(id_{F(M)} \o' F_2(N,P)))(\a^{-i'+1}_M(m) \o' (n \o' \a^{j'-1}_P(p)))\\
&=\a^{i-2i'+1}_M(m) \o (\a_N^{i+j-i'-j'}(n) \o \a^{2j-j'-1}_P(p))\\
&=F(a_{M,N,P})(\a^{2i-2i'}_M(m) \o \a_N^{i+j-i'-j'}(n)) \o \a^{j-j'}_P(p)\\
&=(F(a_{M,N,P})\circ F_2(M \o N, P)\circ(F_2(M.N) \o' id_P))((m\o' n)\o' p).
\endaligned$$

At last, it is easy to get that
$$\aligned
F(l_M)\circ F_2(k,M)\circ (F_0 \o' id_M)=l'_{F(M)},
\endaligned$$
and
$$\aligned
F(r_M)\circ F_2(M,k)\circ (id_M \o' F_0)=r'_{F(M)},
\endaligned$$
hence $F=(F,F_2,F_0)$ is a monoidal functor.

Obviously $F$ is invertible, thus the conclusion holds.

Furthermore, if $H$ is a quasitriangular monoidal Hom-bialgebra with the $R$-matrix $R \in H \o H$, for any $m \in M$, $n \in N$, we compute
\begin{eqnarray*}
&&~~~~(F_2(N,M)\ci c'_{F(M),F(N)})(m \o' n)\\
&&=F_2(N,M)(\a^{i'}(R^{(2)})\c \a_N^{i'-j'-1}(n) \o \a^{j'}(R^{(1)}) \c \a_M^{j'-i'-1}(m))\\
&&=\a^{i}(R^{(2)})\c \a_N^{i-j'-1}(n) \o \a^{j}(R^{(1)}) \c \a_M^{j-i'-1}(m) \\
&&=F(c_{M,N})(\a_M^{i-i'}(m) \o \a_N^{j-j'}(n))=(F(c_{M,N})\ci F_2(M,N))(m \o' n),
\end{eqnarray*}
which implies $F$ is a braided monoidal functor.
$\hfill \Box$

 \section{\textbf{The Drinfeld twists for monoidal Hom-bialgebras}}
\setcounter{section}{4}

 Let $(H,\alpha,m,\eta,\Delta,\varepsilon)$ be a monoidal Hom-bialgebra.

\begin{defn}
A \emph{Drinfeld twist} for $H$ is an invertible element $\sigma \in H \o H$
such that
\\
$\left\{\begin{array}{l}
(T1)~~(\a \o \a)\sigma = \sigma;\\
(T2)~~(\varepsilon \o id_H)(\sigma) = (id_H \o \varepsilon)(\sigma) = 1_H;\\
(T3)~~\sigma^{(1)} \o \bar{\sigma}^{(1)}\sigma^{(2)}_1 \o \bar{\sigma}^{(2)}\sigma^{(2)}_2 = \bar{\sigma}^{(1)}\sigma^{(1)}_1 \o \bar{\sigma}^{(2)}\sigma^{(1)}_2 \o \sigma^{(2)},
\end{array}\right.$
\\
where $\sigma=\sigma^{(1)} \o \sigma^{(2)}=\bar{\sigma}^{(1)} \o \bar{\sigma}^{(2)}$.
\end{defn}

Note that if $\a = id$, then our definition of Drinfeld twists for a bialgebra is inverse with Drinfeld's (see \cite{dr} and \cite{eg}).
\\

{\bf Notation.}
We write $\sigma^{-1} =\varrho= \varrho^{(1)} \o \varrho^{(2)} = \bar{\varrho}^{(1)} \o \bar{\varrho}^{(2)}=\cdots \in H \o H$, and $\sigma_{21}=\sigma^{(2)} \o \sigma^{(1)}$.

\begin{lem}
1). $\varrho$ satisfies
\begin{equation}
(\a \o \a)(\varrho) = \varrho~;
\end{equation}
2). $\varrho$ satisfies the 2-cocycle condition
\begin{equation}
\varrho^{(1)} \o  \varrho^{(2)}_1\bar{\varrho}^{(1)} \o \varrho^{(2)}_2\bar{\varrho}^{(2)} = \varrho^{(1)}_1\bar{\varrho}^{(1)} \o \varrho^{(1)}_2\bar{\varrho}^{(2)} \o \varrho^{(2)}.
\end{equation}
\end{lem}

{\bf Proof.}
1). We compute as follows
\begin{eqnarray*}
(\a \o \a)(\varrho)&=& 1_H\varrho^{(1)} \o 1_H\varrho^{(2)}\\
&=&(\a^{-2}(\bar{\varrho}^{(1)})\a^{-2}(\sigma^{(1)}))\varrho^{(1)} \o (\a^{-2}(\bar{\varrho}^{(2)})\a^{-2}(\sigma^{(2)}))\varrho^{(2)}\\
&\stackrel {(T1)}{=} &\a^{-1}(\bar{\varrho}^{(1)})(\a^{-1}(\sigma^{(1)})\a^{-1}(\varrho^{(1)})) \o \a^{-1}(\bar{\varrho}^{(2)})(\a^{-1}(\sigma^{(2)})\a^{-1}(\varrho^{(2)}))\\
&=&\a^{-1}(\bar{\varrho}^{(1)})1_H \o \a^{-1}(\bar{\varrho}^{(2)})1_H = \varrho.
\end{eqnarray*}

2).
Firstly, we multiply by $\varrho^{(1)} \o  \varrho^{(2)}_1\bar{\varrho}^{(1)} \o \varrho^{(2)}_2\bar{\varrho}^{(2)}$ on the left to Eq.(T3).
Since
$$\aligned
&~~~~(\varrho^{(1)} \o  \varrho^{(2)}_1\bar{\varrho}^{(1)} \o \varrho^{(2)}_2\bar{\varrho}^{(2)})(\sigma^{(1)} \o \bar{\sigma}^{(1)}\sigma^{(2)}_1 \o \bar{\sigma}^{(2)}\sigma^{(2)}_2)\\
&=\varrho^{(1)}\sigma^{(1)} \o \a(\varrho^{(2)}_1)((\a^{-1}(\bar{\varrho}^{(1)})\a^{-1}(\bar{\sigma}^{(1)}))\sigma^{(2)}_1)
 \o \a(\varrho^{(2)}_2)((\a^{-1}(\bar{\varrho}^{(2)}) \a^{-1}(\bar{\sigma}^{(2)}))\sigma^{(2)}_2) \\
&=\varrho^{(1)}\sigma^{(1)} \o \Delta(\a(\varrho^{(2)}\sigma^{(2)}))=1_H \o 1_H \o 1_H,
\endaligned$$
we have
\begin{equation}
1_H \o 1_H \o 1_H = (\varrho^{(1)} \o  \varrho^{(2)}_1\bar{\varrho}^{(1)} \o \varrho^{(2)}_2\bar{\varrho}^{(2)})(\bar{\sigma}^{(1)}\sigma^{(1)}_1 \o \bar{\sigma}^{(2)}\sigma^{(1)}_2 \o \sigma^{(2)}).
\end{equation}

Secondly, we multiply by $\tilde{\varrho}^{(1)}_1\dot{\varrho}^{(1)} \o \tilde{\varrho}^{(1)}_2\dot{\varrho}^{(2)} \o \tilde{\varrho}^{(2)}$ on the right to Eq.(4.3).
We compute
$$\aligned
&~~~~((\varrho^{(1)} \o  \varrho^{(2)}_1\bar{\varrho}^{(1)} \o \varrho^{(2)}_2\bar{\varrho}^{(2)})(\bar{\sigma}^{(1)}\sigma^{(1)}_1 \o \bar{\sigma}^{(2)}\sigma^{(1)}_2 \o \sigma^{(2)}))(\tilde{\varrho}^{(1)}_1\dot{\varrho}^{(1)} \o \tilde{\varrho}^{(1)}_2\dot{\varrho}^{(2)} \o \tilde{\varrho}^{(2)}) \\
&=(\a(\varrho^{(1)}) \o  \a(\varrho^{(2)}_1\bar{\varrho}^{(1)}) \o \a(\varrho^{(2)}_2\bar{\varrho}^{(2)}))((\bar{\sigma}^{(1)}\sigma^{(1)}_1 \o \bar{\sigma}^{(2)}\sigma^{(1)}_2 \o \sigma^{(2)})\\
&~~~~~~~~(\a^{-1}(\tilde{\varrho}^{(1)}_1\dot{\varrho}^{(1)}) \o \a^{-1}(\tilde{\varrho}^{(1)}_2\dot{\varrho}^{(2)}) \o \a^{-1}(\tilde{\varrho}^{(2)})))\\
&=(\varrho^{(1)} \o  \varrho^{(2)}_1\bar{\varrho}^{(1)}) \o \varrho^{(2)}_2\bar{\varrho}^{(2)})(\bar{\sigma}^{(1)}\dot{\varrho}^{(1)} \o \bar{\sigma}^{(2)}\dot{\varrho}^{(2)} \o 1_H)~~\mbox{(since~Eq.(T1)~and~(4.1))}\\
&=\varrho^{(1)} \o  \varrho^{(2)}_1\bar{\varrho}^{(1)} \o \varrho^{(2)}_2\bar{\varrho}^{(2)},
\endaligned$$
thus Eq.(4.2) holds. $\hfill \Box$

\begin{exam}
Let $A$ be a bialgebra over $k$, $\sigma \in A \o A$ be the usual normalized Drinfeld twist for $A$, $\a:A\rightarrow A$ be an invertible bialgebra homomorphism. If $\sigma$ satisfies $(\a \o \a)(\sigma)=\sigma$,
then $\sigma$ is a Drinfeld twist for ${}^\a A$.

Conversely, if $(H,\a)$ is a monoidal Hom-bialgebra endowed with a Drinfeld twist $\sigma$, then $\sigma$ is a normalized Drinfeld twist for $k$-bialgebra ${}_\a H$.
\end{exam}

The following property is a generalization of (Theorem 1.3, \cite{zj}).

\begin{prop}
1). If $(A,\a_A,m_A,\eta_A)$ is an algebra in $Rep^{0,0}(H)$ (the \emph{$H$-Hom-module algebra}), define a new multiplication by
$$a \ci b := (\varrho^{(1)}\c a)(\varrho^{(2)}\c b),$$
for any $a, b \in A$, then $(A,\a_A^2,\ci,1_A)$ is a monoidal Hom-algebra over $k$.

2).If $(C,\a_C,\Delta_C,\varepsilon_C)$ is a coalgebra in $Rep^{0,0}(H)$ (the \emph{$H$-Hom-module coalgebra}), define a new comultiplication by
$$\widehat{\Delta}_C(c) := c_{(1)} \o c_{(2)} = \sigma^{(1)}\c c_1 \o \sigma^{(2)}\c c_2,$$
for any $c \in C$, then $(C,\widehat{\Delta}_C,\varepsilon_C)$ is a coassociative coalgebra over $k$.
\end{prop}

{\bf Proof.}
1).
It is easy to get that $A$ is a left $H$-Hom-module algebra iff $(A,\a_A)$ is both a monoidal Hom-algebra over $k$ and a left $H$-Hom module, such that
\\
$\left\{\begin{array}{l}
(1)~~h \c (ab) = (h_1 \c a)(h_2 \c b);\\
(2)~~h \c 1_A = \varepsilon(h)1_A,\\
\end{array}\right.$
\\
for any $a,b \in A$ and $h \in H$. Thus we have
\begin{eqnarray*}
(a \ci b)\ci \a_A^2(c)&=&( (\varrho^{(1)}_1\c(\bar{\varrho}^{(1)}\c a)) (\varrho^{(1)}_2\c(\bar{\varrho}^{(2)}\c b)) )(\varrho^{(2)}\c \a_A^2(c))\\
&\stackrel {(4.2)}{=}&((\a^{-1}(\varrho^{(1)})\c \a_A(a)) (\a^{-1}(\varrho^{(2)}_1\bar{\varrho}^{(1)})\c\a_A(b)))((\varrho^{(2)}_2\bar{\varrho}^{(2)})\c \a_A^2(c))\\
&=&(\varrho^{(1)}\c \a^2_A(a))( (\varrho^{(2)}_1 \c (\bar{\varrho}^{(1)} \c b)) (\varrho^{(2)}_2 \c (\bar{\varrho}^{(2)} \c c)) )\\
&=& a^2_A(a) \ci(b \ci c),
\end{eqnarray*}
where $a,b,c \in A$, and
$$\aligned
1_A \ci a &= (\varrho^{(1)}\c 1_A)(\varrho^{(2)}\c a)\\
&= 1_A(\varepsilon(\varrho^{(1)})\varrho^{(2)}\c a ) =a^2_A(a)\\
&=a \ci 1_A.
\endaligned$$
So $(A,\a_A^2,\ci,1_A)$ is a monoidal Hom-algebra over $k$.

2). Note that $C$ is a left $H$-Hom-module coalgebra if and only if $(C,\a_C)$ is both a monoidal Hom-coalgebra over $k$ and a left $H$-Hom-module,
such that
\\
$\left\{\begin{array}{l}
(1)~~\Delta_C(h \c c) = h_1 \c c_1 \o h_2 \c c_2;\\
(2)~~\varepsilon_C(h \c c) = \varepsilon(h)\varepsilon_C(c),\\
\end{array}\right.$
\\
for any $c \in C$. Thus we have
\begin{eqnarray*}
\widehat{\Delta}_C(c_{(1)}) \o c_{(2)}&=& \bar{\sigma}^{(1)}\c(\sigma^{(1)}_1\c c_{11}) \o \bar{\sigma}^{(2)}\c(\sigma^{(1)}_2\c c_{12})
\o \sigma^{(2)}\c c_{2} \\
&\stackrel {(T3)}{=}& \sigma^{(1)}\c \a_C(c_{11}) \o (\bar{\sigma}^{(1)}\sigma^{(2)}_1)\c \a_C(c_{12}) \o (\bar{\sigma}^{(2)}\sigma^{(2)}_2)\c c_{2}\\
&=& \sigma^{(1)}\c c_{1} \o (\sigma^{(2)}\c c_{2})_{(1)} \o (\sigma^{(2)}\c c_{2})_{(2)}\\
&=&c_{(1)} \o \widehat{\Delta}_C(c_{(2)}),
\end{eqnarray*}
and
$$\aligned
\varepsilon_C(c_{(1)})c_{(2)}&=\varepsilon(\sigma^{(1)})\sigma^{(2)}\c(\varepsilon_C(c_1)c_2)\\
&=c=c_{(1)}\varepsilon_C(c_{(2)}).
\endaligned$$
Hence $(C,\widehat{\Delta}_C,\varepsilon_C)$ is a coassociative coalgebra over $k$.
$\hfill \Box$
\\

For any $x \in H$, define a new comultiplication $\Delta^\sigma$ on $H$ by
\begin{eqnarray}\Delta^\sigma(x)=x_{[1]} \o x_{[2]} = (\sigma\Delta(x))\varrho.\end{eqnarray}

Note that it is easy to get $\Delta^\sigma(x) = (\sigma\Delta(x))\varrho=\sigma(\Delta(x)\varrho)$.

\begin{lem}\label{hm}
$\Delta^\sigma$ is a Hom-algebra map preserving unit.
\end{lem}

{\bf Proof.}
For one thing, obviously $\Delta^\sigma$ preserves unit and satisfies $\Delta^\sigma(\a(x)) = (\a \o \a)\Delta^\sigma(x)$ for any $x \in H$.

For another, for any $x,y \in H$, we have
$$\aligned
&~~~~x_{[1]}y_{[1]} \o x_{[1]}y_{[1]}\\
&=((\sigma^{(1)}x_1)\varrho^{(1)})((\bar{\sigma}^{(1)}y_1)\bar{\varrho}^{(1)}) \o ((\sigma^{(2)}x_2)\varrho^{(2)})((\bar{\sigma}^{(2)}y_2)\bar{\varrho}^{(2)})\\
&=((\sigma^{(1)}x_1)\a^{-1}(\varrho^{(1)}\bar{\sigma}^{(1)}))\a(y_1\bar{\varrho}^{(1)}) \o ((\sigma^{(2)}x_2)\a^{-1}(\varrho^{(2)}\bar{\sigma}^{(2)}))\a(y_2\bar{\varrho}^{(2)})\\
&=(\a(\sigma^{(1)})(x_1y_1))\a^2(\bar{\varrho}^{(1)}) \o (\a(\sigma^{(2)})(x_2y_2))\a^2(\bar{\varrho}^{(2)}) \\
&=(xy)_{[1]} \o (xy)_{[2]}.~~\mbox{(since~Eq.(T1)~and~(4.1))}
\endaligned$$

Thus the conclusion holds.
$\hfill \Box$

\begin{thm}
$H^\sigma=(H,\a,m,\eta,\Delta^\sigma,\varepsilon)$ is a Hom-bialgebra.
\end{thm}

{\bf Proof.}
Firstly, for any $x \in H$, we have
$$\aligned
x_{[1]}\varepsilon(x_{[2]})&=(\sigma^{(1)}x_1)\varrho^{(1)}\varepsilon(\sigma^{(2)})\varepsilon(x_2)\varepsilon(\varrho^{(2)})\\
&=(1_H\a^{-1}(x))1_H = \a(x)\\
&=\varepsilon(\sigma^{(1)})\varepsilon(x_1)\varepsilon(\varrho^{(1)})(\sigma^{(2)}x_2)\varrho^{(2)}=\varepsilon(x_{[1]})x_{[2]}.
\endaligned$$

Secondly, we have
\begin{eqnarray*}
&&~~~~\Delta^\sigma(x_{[1]}) \o \a(x_{[2]})\\
&=&(\bar{\sigma}^{(1)}((\sigma_1^{(1)}x_{11})\varrho^{(1)}_1))\bar{\varrho}^{(1)} \o (\bar{\sigma}^{(2)}((\sigma_2^{(1)}x_{12})\varrho^{(1)}_2))\bar{\varrho}^{(2)} \o \a((\sigma^{(2)}x_2)\varrho^{(2)})\\
&=&((\a^{-1}(\bar{\sigma}^{(1)})\a(\sigma_1^{(1)}))\a^2(x_{11}))(\a(\varrho^{(1)}_1)\a^{-1}(\bar{\varrho}^{(1)})) \\
&&~~~~~~~~\o ((\a^{-1}(\bar{\sigma}^{(2)})\a(\sigma_2^{(1)}))\a^2(x_{12}))(\a(\varrho^{(1)}_2)\a^{-1}(\bar{\varrho}^{(2)})) \o \a((\sigma^{(2)}x_2)\varrho^{(2)})\\
&=&((\bar{\sigma}^{(1)}\sigma_1^{(1)})\a^2(x_{11}))(\varrho^{(1)}_1\bar{\varrho}^{(1)}) \o ((\bar{\sigma}^{(2)}\sigma_2^{(1)})\a^2(x_{12}))(\varrho^{(1)}_2\bar{\varrho}^{(2)}) \o (\sigma^{(2)}\a (x_2))\varrho^{(2)}\\
&\stackrel {(T3),(4.2)}{=}&(\sigma^{(1)}\a(x_{1}))\varrho^{(1)} \o ((\bar{\sigma}^{(1)}\sigma^{(2)}_1)\a^2(x_{21}))(\varrho^{(2)}_1\bar{\varrho}^{(1)}) \o ((\bar{\sigma}^{(2)}\sigma^{(2)}_2)\a^2(x_{22}))(\varrho^{(2)}_2\bar{\varrho}^{(2)})\\
&=&\a((\sigma^{(1)}x_1)\varrho^{(1)}) \o ((\a^{-1}(\bar{\sigma}^{(1)})\a(\sigma^{(2)}_1))\a^2(x_{21}))(\a(\varrho^{(2)}_1)\a^{-1}(\bar{\varrho}^{(1)}))  \\
&&~~~~~~~~\o ((\a^{-1}(\bar{\sigma}^{(2)})\a(\sigma^{(2)}_2))\a^2(x_{22}))(\a(\varrho^{(2)}_2)\a^{-1}(\bar{\varrho}^{(2)})) \\
&=&\a((\sigma^{(1)}x_1)\varrho^{(1)}) \o (\bar{\sigma}^{(1)}((\sigma^{(2)}_1x_{21})\varrho^{(2)}_1))\bar{\varrho}^{(1)} \o (\bar{\sigma}^{(2)}((\sigma^{(2)}_2x_{22})\varrho^{(2)}_2))\bar{\varrho}^{(2)}\\
&=&\a(x_{[1]}) \o\Delta^\sigma(x_{[2]}).
\end{eqnarray*}

Finally, we have
\begin{eqnarray*}
\Delta^\sigma(\a(x))&=&(\sigma^{(1)}\a(x_{1}))\varrho^{(1)} \o (\sigma^{(2)}\a(x_{2}))\varrho^{(2)} \\
&=&\a((\sigma^{(1)}x_{1})\varrho^{(1)}) \o \a((\sigma^{(2)}x_{2})\varrho^{(2)}),~~\mbox{(since~Eq.(T1)~and~(4.1))}
\end{eqnarray*}
which implies $(H, \a, \Delta^\sigma,\varepsilon)$ is a Hom-coalgebra.

Since $(H,\a,m,\eta)$ is already a Hom-algebra, by Lemma \ref{hm}, the conclusion holds.
$\hfill \Box$
\\

For any given monoidal Hom-bialgebra $(H,\a)$ endowed with a Drinfeld twist $\sigma \in H \o H$, we already know that ${}_\a H$ is a bialgebra and $\sigma$ is a Drinfeld twist on it. Thus $({}_\a H)^\sigma$ is a new bialgebra which coproduct is also given by Eq.(4.4).
By Example \ref{bh}, $((({}_\a H)^\sigma)^\a,\a)$ is a Hom-bialgebra.

\begin{thm}
$(({}_\a H)^\sigma)^\a = H^\sigma$ as Hom-bialgebras. Furthermore, $({}_\a H)^\a=H^\sigma$ if and only if $\sigma=1_H \o 1_H$.
\end{thm}

{\bf Proof.}
We only need to show $\Delta_{(({}_\a H)^\sigma)^\a} = \Delta_{H^\sigma}$.

In fact, denote the multiplication $m_{{}_\a H}(x \o y) = \a^{-1}(x)\a^{-1}(y)$ by $x \ast y$ for any $x,y \in H$. We have
$$\aligned
\Delta_{(({}_\a H)^\sigma)^\a}(x)&=\Delta_{({}_\a H)^\sigma}(\a(x))\\
&=(\sigma \ast \Delta_{{}_\a H}(\a(x))) \ast \varrho\\
&=(\sigma \ast \Delta(\a^2(x))) \ast \varrho \\
&=(\sigma \Delta(x)) \varrho.
\endaligned$$
Thus the conclusion holds.
$\hfill \Box$

\begin{thm}
If $(H,\a,S)$ is a monoidal Hom-Hopf algebra,
then $H^\sigma$ is a Hom-Hopf algebra.
\end{thm}

{\bf Proof.}
Firstly, since
$$\aligned
&(1_H \o \bar{\varrho}^{(1)}\o \bar{\varrho}^{(2)}) ( (\sigma^{(1)} \o \bar{\sigma}^{(1)}\sigma^{(2)}_1 \o \bar{\sigma}^{(2)}\sigma^{(2)}_2 )(\a(\varrho^{(1)}_1) \o \a(\varrho^{(1)}_2) \o \varrho^{(2)})) \\
&~~~~~~~~= (1_H \o \bar{\varrho}^{(1)}\o \bar{\varrho}^{(2)}) ( (\bar{\sigma}^{(1)}\sigma^{(1)}_1 \o \bar{\sigma}^{(2)}\sigma^{(1)}_2 \o \sigma^{(2)}) (\a(\varrho^{(1)}_1) \o \a(\varrho^{(1)}_2) \o \varrho^{(2)}) ),
\endaligned$$
we have
\begin{eqnarray}
\a(\sigma^{(1)}) \o \varrho^{(1)} \sigma^{(2)} \o \a(\varrho^{(2)}) = \sigma^{(1)}\a(\varrho^{(1)}_1) \o \a(\sigma^{(2)}_1)\a(\varrho^{(1)}_2) \o \a(\sigma^{(2)}_2)\varrho^{(2)}.
\end{eqnarray}

Define $S^\sigma: H\rightarrow H$ by
$$
S^\sigma(x) = (\sigma^{(1)} (S(\a^{-1}(\sigma^{(2)})) (S(\a^{-4}(x)) S(\a^{-3}(\varrho^{(1)}))) ) )\varrho^{(2)},
$$
for any $x \in H$, where $p,q \in \mathbb{Z}$. Obviously $S \ci \a = \a \ci S$.

Since
\begin{eqnarray*}
&&~~~~S^\sigma(x_{[1]})x_{[2]} \\
&=&((\sigma^{(1)} (S(\a^{-1}(\sigma^{(2)})) ((S(\a^{-4}(\varrho^{(1)}))(S(\a^{-4}(x_1))S(\a^{-4}(\sigma^{(1)})))) \\
&&~~~~~~~~S(\a^{-3}(\varrho^{(1)})) ) ) )\varrho^{(2)}) ((\sigma^{(2)} x_2)\varrho^{(2)})\\
&=&( (\sigma^{(1)} ( S(\a^{-1}(\sigma^{(2)})) ( (S(\a^{-4}(\varrho^{(1)})) S(\a^{-3}(x_1)) )S(\a^{-4}(\varrho^{(1)})\a^{-3}(\sigma^{(1)}))  ) )) \\
&&~~~~~~~~( \a^{-1}(\varrho^{(2)}) \sigma^{(2)}) ) (\a(x_2)\varrho^{(2)})\\
&=& ( \a(\sigma^{(1)}) ( S(\sigma^{(2)}) ( (S(\a^{-2}(\varrho^{(1)})) S(\a^{-1}(x_1)) )  ) )) (\a(x_2)\varrho^{(2)})\\
&=& ( ( \sigma^{(1)} ( S(\a^{-1}(\sigma^{(2)})) S(\a^{-2}(\varrho^{(1)})) ) ) S(\a(x_1)) ) (\a(x_2)\varrho^{(2)}) \\
&=& ( ( \sigma^{(1)} ( S(\a^{-1}(\sigma^{(2)})) S(\a^{-2}(\varrho^{(1)})) ) )  (S(x_1) x_2) ) \a(\varrho^{(2)}) \\
&=&( \a(\sigma^{(1)}) ( S(\sigma^{(2)}) S(\a^{-1}(\varrho^{(1)})) ) ) \a(\varrho^{(2)})\varepsilon(x) \\
&=&( \a(\sigma^{(1)}) ( S(\sigma^{(2)}) S(\varrho^{(1)}) ) ) \a^{2}(\varrho^{(2)})\varepsilon(x) \\
&=&( \a(\sigma^{(1)}) S(\varrho^{(1)}\sigma^{(2)}) ) \a^{2}(\varrho^{(2)})\varepsilon(x) \\
&\stackrel {(4.5)}{=}&( ( \sigma^{(1)} \a^{1}(\varrho^{(1)}_1) )  ( S(\a(\varrho^{(1)}_2)) S(\a(\sigma^{(2)}_1)) ) ) ( \a^{2}(\sigma^{(2)}_2) \a(\varrho^{(2)}) )\varepsilon(x) \\
&=&( ( \sigma^{(1)} ( \varrho^{(1)}_1 S(\a(\varrho^{(1)}_2)) ) ) S(\a^{2}(\sigma^{(2)}_1)) ) ( \a^{2}(\sigma^{(2)}_2) \a(\varrho^{(2)}) )\varepsilon(x) \\
&=& \a^{2}(\sigma^{(1)})( S(\a^{2}(\sigma^{(2)}_1)) \a^{2}(\sigma^{(2)}_2) )\varepsilon(x) \\
&=& \a^{3}(\sigma^{(1)})\varepsilon(\sigma^{(2)})\varepsilon(x) = 1_H\varepsilon(x),
\end{eqnarray*}
which implies $S^\sigma$ is the antipode of $H^\sigma$.
$\hfill \Box$

\begin{prop}\label{mhb}
If $(H,\a)$ is a quasitiangular monoidal Hom-bialgebra with the $R$-matrix $R = R^{(1)} \o R^{(1)}$,
then $H^\sigma$ is a quasitriangular Hom-bialgebra.
\end{prop}

{\bf Proof.} Define $R^\sigma = (\sigma_{21}R)\varrho \in H \o H$, we will prove $R^\sigma$ is an $R$-matrix in $H^\sigma$.

Firstly, we will check Eq.(q2). For any $x \in H$, we have
\begin{eqnarray*}
&&~~~~{\Delta^\sigma}^{op}(x)R^\sigma\\
&=&((\sigma^{(2)}x_2)\varrho^{(2)})((\bar{\sigma}^{(2)}R^{(1)})\bar{\varrho}^{(1)}) \o ((\sigma^{(1)}x_1)\varrho^{(1)})((\bar{\sigma}^{(1)}R^{(2)})\bar{\varrho}^{(2)})\\
&=&(\sigma^{(2)}(x_2((\varrho^{(2)}\bar{\sigma}^{(2)})R^{(1)})))\bar{\varrho}^{(1)} \o (\sigma^{(1)}(x_1((\varrho^{(1)}\bar{\sigma}^{(1)})R^{(2)})))\bar{\varrho}^{(2)} \\
&\stackrel {(Q2)}{=}&(\sigma^{(2)}(R^{(1)}x_1))\bar{\varrho}^{(1)} \o (\sigma^{(1)}(R^{(2)}x_2))\bar{\varrho}^{(2)}\\
&=&((\sigma^{(2)}R^{(1)})\varrho^{(1)})((\bar{\sigma}^{(1)}x_1)\bar{\varrho}^{(1)}) \o ((\sigma^{(1)}R^{(2)})\varrho^{(2)})((\bar{\sigma}^{(2)}x_2)\bar{\varrho}^{(2)})\\
&=&R^\sigma\Delta^\sigma(x).
\end{eqnarray*}

Secondly, to verify Eq.(q3), we compute
\begin{eqnarray*}
&&~~~~\Delta^\sigma({R^\sigma}^{(1)}) \o \a({R^\sigma}^{(2)}) \\
&=&(\bar{\sigma}^{(1)}((\sigma^{(2)}_1 R^{(1)}_1)\varrho^{(1)}_1))\bar{\varrho}^{(1)} \o (\bar{\sigma}^{(2)}((\sigma^{(2)}_2 R^{(1)}_2)\varrho^{(1)}_2))\bar{\varrho}^{(2)} \o \a((\sigma^{(1)} R^{(2)})\varrho^{(2)})\\
&\stackrel {(Q3)}{=}&((\bar{\sigma}^{(1)}\sigma^{(2)}_1)\a^2(R^{(1)}))(\varrho^{(1)}_1\bar{\varrho}^{(1)}) \o ((\bar{\sigma}^{(2)}\sigma^{(2)}_2)\a^2(r^{(1)}))(\varrho^{(1)}_2\bar{\varrho}^{(2)}) \o (\sigma^{(1)}\\
&&~~~~~~~~~~(\a(R^{(2)}r^{(2)})))\varrho^{(2)}   \\
&\stackrel {(T3),(4.2)}{=}& ( \a(\bar{\sigma}^{(2)})(\sigma^{(1)}_2 R^{(1)}) )\varrho^{(1)} \o (\sigma^{(2)} (r^{(1)} \a^{-1}(\varrho^{(2)}_1) ) ) \a(\bar{\varrho}^{(1)}) \o (\a(\bar{\sigma}^{(1)}) (\sigma^{(1)}_1 R^{(2)}) ) \\
&&~~~~~~~~~~( (r^{(2)}\a^{-1}(\varrho^{(2)}_2) ) \bar{\varrho}^{(2)})\\
&\stackrel {(Q2)}{=}&(\bar{\sigma}^{(2)}(R^{(1)}\sigma^{(1)}_1))\a(\varrho^{(1)}) \o (\sigma^{(2)}(\varrho^{(2)}_2 r^{(1)}))\a(\bar{\varrho}^{(1)}) \o(\bar{\sigma}^{(1)}(R^{(2)}\sigma^{(1)}_2))\\
&&~~~~~~~~((\varrho^{(2)}_1 r^{(2)})\bar{\varrho}^{(2)})\\
&=&\a((\a^{-1}(\bar{\sigma}^{(2)})R^{(1)})(\sigma^{(1)}_1\a^{-1}(\varrho^{(1)}))) \o (\sigma^{(2)} \a(\varrho^{(2)}_2))(\a(r^{(1)}) \bar{\varrho}^{(1)})\\
&&~~~~~~~~\o ( (\a^{-1}(\bar{\sigma}^{(1)})R^{(2)}) (\sigma^{(1)}_2 \varrho^{(2)}_1) )( \a(r^{(2)}) \bar{\varrho}^{(2)} ) \\
&=&\a((\a^{-1}(\bar{\sigma}^{(2)})R^{(1)})((\dot{\varrho}^{(1)}\dot{\sigma}^{(1)})\a^{-1}(\sigma^{(1)}_1\a^{-1}(\varrho^{(1)})))) \o (\sigma^{(2)} \a(\varrho^{(2)}_2))(\a(r^{(1)}) \bar{\varrho}^{(1)})\\
&&~~~~~~~~\o ( (\a^{-1}(\bar{\sigma}^{(1)})R^{(2)}) ((\dot{\varrho}^{(2)}\dot{\sigma}^{(2)}) \a^{-1}(\sigma^{(1)}_2 \varrho^{(2)}_1) ))( \a(r^{(2)}) \bar{\varrho}^{(2)} ) \\
&\stackrel {(T3)}{=}&\a((\a^{-1}(\bar{\sigma}^{(2)})R^{(1)})((\dot{\varrho}^{(1)}\a^{-2}(\sigma^{(1)}))\a^{-1}(\varrho^{(1)}))) \o ((\dot{\sigma}^{(2)}\sigma^{(2)}_2)\a(\varrho^{(2)}_2))(\a(r^{(1)}) \bar{\varrho}^{(1)}) \\
&&~~~~~~~~ \o ( (\a^{-1}(\bar{\sigma}^{(1)})R^{(2)}) ( (\dot{\varrho}^{(2)}  \a^{-2}(\dot{\sigma}^{(1)}\sigma^{(2)}_1) ) \varrho^{(2)}_1)  ) ( \a(r^{(2)}) \bar{\varrho}^{(2)}) \\
&=&\a((\bar{\sigma}^{(2)}R^{(1)})\a(\dot{\varrho}^{(1)})) \o \a^2(\dot{\sigma}^{(2)})(r^{(1)}\bar{\varrho}^{(1)}) \o ( (\bar{\sigma}^{(1)}R^{(2)}) (\dot{\varrho}^{(2)}\dot{\sigma}^{(1)}) )(r^{(2)}\bar{\varrho}^{(2)}) \\
&=&(\a(\bar{\sigma}^{(2)})\a(R^{(1)}))\a^2(\dot{\varrho}^{(1)}) \o (\a(\dot{\sigma}^{(2)})r^{(1)})\a(\bar{\varrho}^{(1)}) \\
&&~~~~~~~~\o \a^{-1}( ((\a(\bar{\sigma}^{(1)})\a(R^{(2)}))\a^2(\dot{\varrho}^{(2)})) ((\a(\dot{\sigma}^{(1)})r^{(2)})\a(\bar{\varrho}^{(2)})))\\
&=&\a({R^\sigma}^{(1)}) \o \a({r^\sigma}^{(1)}) \o {R^\sigma}^{(2)}{r^\sigma}^{(2)}.
\end{eqnarray*}

Eq.(q4) could be obtained in a similar way.

At last, since $(\a \o\a)(R^\sigma) = R^\sigma$, $(H^\sigma,\a,R^\sigma)$ is a quasitriangular Hom-bialgebra.
$\hfill \Box$
\\

Recall from \cite{zxh} that if $(B,\a_B)$ is a Hom-bialgebra over $k$ and $\a_B$ is invertible, then
the representation category $Rep^{i,j}(B)$ of $B$ is a monoidal category with the following structure:

$\bullet$ the objects are left $B$-Hom-modules and the morphism $f:(M,\a_M)\rightarrow(N,\a_N)$ is $B$-linear and satisfies $f \ci \a_M = \a_N \ci f$;

$\bullet$ the tensor product of $(M,\a_M)$ and $(N,\a_N)$ is $(M \o N, \a_M \o \a_N)$, with the $B$-action given by
$$
b \cdot (m \o n) = \a_B^i(b_1) \cdot m\o \a_B^j(b_2) \cdot n, ~~~~b \in B,~ m\in M,~n\in N;
$$

$\bullet$ the unit object is $(k,id_k)$ with the $B$-action $b\cdot \lambda = \varepsilon(b)\lambda$, where $b \in B$, $\lambda \in k$;

$\bullet$ the associativity constraint $a$ is given by
$$
a_{M,N,P}((m \o n) \o p) = \a_M^{-i-1}(m) \o (n \o \a_P^{j+1}(p));
$$

$\bullet $ the left unit constraint $l$ is given by
$$
l_M(\lambda \o m) = \lambda \a^{-j-1}_M(m),
$$
for any $m \in M$ and $\lambda \in k$;

$\bullet $ the right unit constraint $r$ is given by
$$
r_M(m \o\lambda) = \lambda \a^{-i-1}_M(m),
$$
for any $m \in M$ and $\lambda \in k$.

Furthermore, if $(B,\a_B)$ is a quasitriangular Hom-bialgebra with the $R$-matrix $R=R^{(1)} \o R^{(2)}$, then $Rep^{i,j}(B)$ is a braided category with the braiding
$$
c_{M,N}(m \o n) = \a_B^i(R^{(2)})\cdot \a_N^{i-j-1}(n) \o \a_B^j(R^{(1)})\cdot \a_M^{j-i-1}(m),
$$
for any $m \in M$, $n \in N$.

\begin{thm}
$Rep^{i+3,j+3}(H)$ and $Rep^{i,j}(H^\sigma)$ are isomorphic as monoidal categories. Moreover, if $H$ is a quasitriangular monoidal Hom-bialgebra, then $Rep^{i+3,j+3}(H)$ and $Rep^{i,j}(H^\sigma)$ are braided isomorphic.
\end{thm}

{\bf Proof.}
Define a functor
$$G=(G,G_2,G_0):(Rep^{i+3,j+3}(H),\o, k,a,l,r)\rightarrow (Rep^{i,j}(H^\sigma),\bar{\o},k,\bar{a},\bar{l},\bar{r})$$
by
$$G(M):=M,~~G(f):=f,~~G_0=id_k,$$
and
$$G_2(M,N):G(M) \bar{\o} G(N)\rightarrow G(M \o N),~~m \bar{\o} n\mapsto \a^i(\varrho^{(1)}) \c m \o \a^j(\varrho^{(2)}) \c n,$$
where the $H^\sigma$-Hom-module structure of $M$ is given by
$$h\rightharpoonup m = h\c m,~~~h \in H^\sigma,~~~m \in M,$$
for any $m \in M$, $n \in N$, and $f:M\rightarrow N  \in Mor(Rep^{i+3,j+3}(H))$.
Obviously $G_2$ is natural and compatible with the Hom structure map.

Firstly, since $G_2(M,N)$ is invertible, and
$$\aligned
&~~~~G_2(M,N)(h\rightharpoonup(m \bar{\o} n))\\
&=G_2(M,N)(\a^{i}(h_{[1]})\rightharpoonup m \bar{\o} \a^{j}(h_{[2]})\rightharpoonup n)\\
&=\a^{i+3}(\varrho^{(1)}) \c (\a^{i}((\sigma^{(1)}h_1)\bar{\varrho}^{(1)}) \c m) \o \a^{j+3}(\varrho^{(2)}) \c (\a^{j}((\sigma^{(2)}h_2)\bar{\varrho}^{(2)})\c n)\\
&=((\a^{i+1}(\varrho^{(1)})\a^{i+1}(\sigma^{(1)}))(\a^{i+1}(h_1)\a^{i}(\bar{\varrho}^{(1)})))\c\a_M(m) \\
&~~~~~~~~\o ((\a^{j+1}(\varrho^{(2)})\a^{j+1}(\sigma^{(2)}))(\a^{j+1}(h_1)\a^{j}(\bar{\varrho}^{(2)})))\c\a_N(n)\\
&=\a^{i+3}(h_1) \c (\a^{i+1}(\bar{\varrho}^{(1)}) \c m) \o \a^{j+3}(h_2) \c (\a^{j+1}(\bar{\varrho}^{(2)}) \c n) \\
&=h\c(\a^i(\bar{\varrho}^{(1)}) \c m \o \a^j(\bar{\varrho}^{(2)}) \c n ))\\
&=h\rightharpoonup(G_2(M,N)( m \bar{\o} n)),
\endaligned$$
$G_2$ is a natural isomorphism in $Rep^{i,j}(H^\sigma)$.

Secondly, we have
\begin{eqnarray*}
&&~~~~(G_2(M,N \o P)\circ(id_{G(M)} \bar{\o} G_2(N,P))\circ(\bar{a}_{G(M),G(N),G(P)}))((m\bar{\o} n)\bar{\o} p)\\
&=&(G_2(M,N \o P)\circ(id_{G(M)} \bar{\o} G_2(N,P)))(\a^{-i-1}_M(m) \bar{\o} (n \bar{\o} \a^{j+1}_P(p)))\\
&=&\a^{i}(\bar{\varrho}^{(1)})\c\a^{-i-1}_M(m) \o \a^{j}(\bar{\varrho}^{(2)})\c(\a^{i}(\varrho^{(1)})\c n \o \a^{j}(\varrho^{(2)}) \c\a^{j+1}_P(p))\\
&=&\a^{i}(\bar{\varrho}^{(1)})\c\a^{-i-1}_M(m) \o ((\a^{i+j+2}(\bar{\varrho}^{(2)}_1)\a^{i}(\varrho^{(1)}))\c \a_N(n) \o (\a^{2j+2}(\bar{\varrho}^{(2)}_2)\a^{j}(\varrho^{(2)})) \c \a^{j+2}_P(p))\\
&=&\a^{i-j-2}(\bar{\varrho}^{(1)})\c\a^{-i-1}_M(m) \o ((\a^{i}(\bar{\varrho}^{(2)}_1)\a^{i}(\varrho^{(1)}))\c \a_N(n) \o (\a^{j}(\bar{\varrho}^{(2)}_2)\a^{j}(\varrho^{(2)})) \c \a^{j+2}_P(p))\\
&\stackrel {(4.2)}{=}&\a^{i-j-2}(\bar{\varrho}^{(1)}_1\varrho^{(1)})\c\a^{-i-1}_M(m) \o (\a^{i}(\bar{\varrho}^{(1)}_2\varrho^{(2)})\c \a_N(n) \o \a^{j}(\bar{\varrho}^{(2)})\c \a^{j+2}_P(p))\\
&=&(\a^i(\bar{\varrho}^{(1)}_1)\a^{-2}(\varrho^{(1)}))\c\a^{-i-1}_M(m) \o ((\a^{i+j+2}(\bar{\varrho}^{(1)}_2)\a^{j}(\varrho^{(2)}))\c \a_N(n) \o (a^{2j+2}(\bar{\varrho}^{(2)}))\c \a^{j+2}_P(p))\\
&=&G(a_{M,N,P})(((\a^{2i+2}(\bar{\varrho}^{(1)}_1)\a^i(\varrho^{(1)})) \c \a_M(m) \o (\a^{i+j+2}(\bar{\varrho}^{(1)}_2)\a^{j}(\varrho^{(2)})) \c \a_N(n)) \o a^{j}(\bar{\varrho}^{(2)})\c p)\\
&=&G(a_{M,N,P})(\a^{i}(\bar{\varrho}^{(1)})\c(\a^{i}(\varrho^{(1)})\c m \o \a^{j}(\varrho^{(2)}) \c n)\o a^{j}(\bar{\varrho}^{(2)})\c p)\\
&=&(G(a_{M,N,P}) \ci G_2(M \o N, P) \ci (G_2(M,N) \bar{\o} id_P))((m\bar{\o} n)\bar{\o} p).
\end{eqnarray*}

At last, it is easy to get
$$\aligned
G(l_M)\circ G_2(k,M)\circ (G_0 \bar{\o} id_M)=\bar{l}_{G(M)},
\endaligned$$
and
$$\aligned
G(r_M)\circ G_2(M,k)\circ (id_M \bar{\o} G_0)=\bar{r}_{G(M)},
\endaligned$$
hence $G=(G,G_2,G_0)$ is a monoidal functor.

Obviously $G$ is invertible, thus $Rep^{i+3,j+3}(H)$ and $Rep^{i,j}(H^\sigma)$ are monoidal isomorphic.

Furthermore, if $R \in H \o H$ is an $R$-matrix, from Proposition \ref{mhb}, $(H^\sigma,\a,R^\sigma)$ is a quasitriangular Hom-bialgebra. Suppose that the braiding in $Rep^{i+3,j+3}(H)$ is $c$ and the braiding in $Rep^{i,j}(H^\sigma)$ is $\bar{c}$, then for any
$m \in M$, $n \in N$, we have
$$\aligned
&~~~~(G_2(N,M)\ci \bar{c}_{G(M),G(N)})(m \bar{\o} n)\\
&=G_2(N,M)(\a^{i}(\sigma^{(1)}(R^{(2)}\varrho^{(2)}))\rightharpoonup \a_N^{i-j-1}(n) \bar{\o} \a^{j}(\sigma^{(2)}(R^{(1)}\varrho^{(1)})) \rightharpoonup \a_M^{j-i-1}(m))\\
&=\a^{i}( (\bar{\varrho}^{(1)} \sigma^{(1)}) (R^{(2)}\varrho^{(2)}) )\c \a_N^{i-j}(n) \o \a^{j}( (\bar{\varrho}^{(2)} \sigma^{(2)})  (R^{(1)}\varrho^{(1)}) ) \c \a_M^{j-i}(m) \\
&=G(c_{M,N})(\a^i(\varrho^{(1)}) \c m \o \a^j(\varrho^{(2)}) \c n)=(G(c_{M,N})\ci G_2(M,N))(m \o' n),
\endaligned$$
which implies $G$ is a braided monoidal functor.
$\hfill \Box$
\\

\begin{center}
 {\bf ACKNOWLEDGEMENT}
\end{center}
The work was partially supported by  the Fundamental Research Funds for the Central Universities  (NO. 3207013906),
and the NSF of China (NO. 11371088), and the NSF of Jiangsu Province (NO. BK2012736).
\\

 \end{document}